\documentclass{elsart}
 \journal{Journal of Combinatorial Theory Series A}

 \usepackage{amsmath}
 \usepackage{amssymb}
 
\numberwithin{equation}{section}

\newcommand{\Bone}{\mathcal{B}_1}
\newcommand{\Gone}{\mathcal{G}_1} 
\newcommand{\Gtwo}{\mathcal{G}_2}
\newcommand{\Gthree}{\mathcal{G}_3}
\newcommand{\Aone}{\mathcal{A}_1}
\newcommand{\Atwo}{\mathcal{A}_2}
\newcommand{\Athree}{\mathcal{A}_3}
\newcommand{\card}{\mathrm{Card}}

 \numberwithin{thm}{section}
\begin{document}
\begin{frontmatter}
 
 \title{A Partition Bijection Related to the \\Rogers-Selberg Identities and \\
 Gordon's Theorem}
 
 \author{Andrew V. Sills}
\address{ 
Department of Mathematical Sciences, Georgia Southern University,
Statesboro, GA 30460}
\ead{asills@GeorgiaSouthern.edu}
\date{March 22, 2007}

 \begin{abstract}
We provide a bijective map from the partitions enumerated by the series side of the
Rogers-Selberg mod 7 identities onto 
partitions associated with a special case of Basil Gordon's combinatorial
generalization of the Rogers-Ramanujan identities.   The implications of 
applying the same map to a special case of David Bressoud's even modulus
analog of Gordon's theorem are also explored.
\end{abstract}

\end{frontmatter}
 
 \section{Introduction}
The celebrated Rogers-Ramanujan identities were first discovered by L.J. Rogers
in 1894~\cite{ljr1894}, but received
little attention until they were
rediscovered independently by S. Ramanujan~\cite[Ch. III, pp. 33 ff.]{pam} and I. Schur~\cite{is}
some two decades later.
They are as follows:
\begin{thm}[The Rogers-Ramanujan Identities]
\begin{equation}\label{rr2}
\sum_{j=0}^\infty \frac{ q^{j^2+j} }{(q;q)_j} = \underset{j\not\equiv 0,\pm 1{\hskip-3mm}\pmod 5}{\prod_{j\geqq 1}}
 \frac{1}{1-q^j}
\end{equation}  and
\begin{equation}\label{rr1}
\sum_{j=0}^\infty \frac{ q^{j^2} }{(q;q)_j} = \underset{j\not\equiv 0,\pm 2{\hskip-3mm}\pmod 5}{\prod_{j\geqq 1}}
 \frac{1}{1-q^j},
\end{equation}
where $(a;q)_j:= \prod_{h=0}^{j-1} (1-aq^h)$ for $j$ a positive integer and $(a;q)_0:=1$.
 \end{thm}
 The Rogers-Ramanujan identities and the other $q$-series identities mentioned in this
paper may be
considered identities of analytic functions that are valid if and only if $|q|<1$.  Since our
interest here is combinatorial, convergence conditions will not be mentioned again.
 
Rogers presented a large number of $q$-series identities 
which resemble~\eqref{rr2} and \eqref{rr1}
in~\cite{ljr1894} and~\cite{ljr1917}.  Among them were a set of three identities associated with
the modulus 7, which received little attention until they were independently rediscovered by
A.~Selberg~\cite{as} and then re-proved by F.J.~Dyson in~\cite{fjd}.
\begin{thm}[The Rogers-Selberg Identities] 
\begin{equation} \label{rs1}
 \sum_{j=0}^\infty \frac{ q^{2j^2+2j}  (-q^{2j+2};q)_\infty }{(q^2;q^2)_j}
 = \underset{j\not\equiv 0,\pm 1{\hskip-3mm}\pmod 7}{\prod_{j\geqq 1}}
 \frac{1}{1-q^j},
\end{equation}
\begin{equation} \label{rs2}
 \sum_{j=0}^\infty \frac{ q^{2j^2+2j}  (-q^{2j+1};q)_\infty }{(q^2;q^2)_j}
 = \underset{j\not\equiv 0,\pm 2{\hskip-3mm}\pmod 7}{\prod_{j\geqq 1}}
 \frac{1}{1-q^j},
\end{equation} and
\begin{equation} \label{rs3}
 \sum_{j=0}^\infty \frac{ q^{2j^2}  (-q^{2j+1};q)_\infty }{(q^2;q^2)_j}
 = \underset{j\not\equiv 0,\pm 3{\hskip-3mm}\pmod 7}{\prod_{j\geqq 1}}
 \frac{1}{1-q^j},
\end{equation}
where $(a;q)_\infty = \prod_{h=0}^\infty (1-aq^h)$.
\end{thm}

Both MacMahon~\cite[Ch. III, pp. 33 ff.]{pam} and Schur~\cite{is} showed that~\eqref{rr2} and \eqref{rr1}
could be interpreted as identities in the theory of integer partitions.  (See \S\ref{Defs} for the
definition of partition and related terms.)
\begin{thm}[The Rogers-Ramanujan Identities--Combinatorial Version]
For $i=1,2$, the number of partitions of $n$ into parts which are nonconsecutive integers
greater than $i-1$ and in which no part is repeated equals the number of partitions of $n$
into parts $\not\equiv 0,\pm i \pmod{5}$.
\end{thm}

In 1961, Basil Gordon~\cite{bg} gave the following generalization of the combinatorial
generalization of the Rogers-Ramanujan 
identities:
\begin{thm}[Gordon's theorem] \label{gordon}
Let $G_{k,i}(n)$ denote the number of partitions of $n$ into parts such that $1$ appears
as a part at most $i-1$ times and the total number of appearances of any two consecutive
integers is at most $k-1$.  Let $C_{k,i}(n)$ denote the number of partitions of $n$
into parts $\not\equiv 0,\pm i \pmod{2k+1}$.  Then $G_{k,i}(n) = C_{k,i}(n)$ for 
$1\leqq i \leqq k$ and all integers $n$.
\end{thm}
\noindent Indeed, it is an elementary exercise to see that \eqref{rr2} and \eqref{rr1} correspond
to the $i=1$ and $i=2$ cases, respectively, of the $k=2$ case of Gordon's theorem
(see, e.g., \cite[p. 290 ff.]{hw}).

 The right hand sides of~\eqref{rs1}, \eqref{rs2}, and \eqref{rs3} are
clearly the generating functions for the partitions enumerated by $C_{3,1}(n)$,
$C_{3,2}(n)$, and $C_{3,3}(n)$ respectively.  Nonetheless, relating the 
partitions enumerated by the $G_{3,i}(n)$
to the left hand sides of~\eqref{rs1}, \eqref{rs2}, and \eqref{rs3} is \emph{not}
a straightforward matter.

In a recent paper, George Andrews~\cite{gea} provided the following partition theoretic
interpretation of~\eqref{rs2}:
\begin{thm}[Andrews]\label{andrews}
Let $A_2(n)$ denote the number of partitions of $n$ such that if
$2j$ is the largest repeated even part, then all positive even integers less than $2j$
also appear at least twice, no odd part less than $2j$ appears, and
no part greater than $2j$ is repeated.  Then $A_2(n) = C_{3,2}(n)$ for
all $n$.
\end{thm}
\begin{pf}
Note that
\[ \frac{q^{2j^2+2j}(-q^{2j+1};q)_\infty}{(q^2;q^2)} 
 = \frac{ q^{2+2+4+4+\cdots+2j+2j} }{(q^2;q^2)} \times (-q^{2j+1};q)_\infty.
\]
By the methods of Euler (cf.~\cite[p. 4 ff.]{top}), 
\[ \frac{ q^{2+2+4+4+\cdots+2j+2j} }{(q^2;q^2)} \]
is the generating function for partitions into $2$'s, $4$'s, $6$'s, $\dots, 2j$'s with
each part appearing at least twice, and $(-q^{2j+1};q)_\infty$ is the generating
function for partitions into distinct parts with each part at least $2j+1$.  Thus, by summing over
all nonnegative $j$, it follows that the left hand side of~\eqref{rs2} is the
generating function for $A_2(n)$.  Again by Euler's method, it is 
immediate that
the right hand side of~\eqref{rs2} is the generating function for $C_{3,2}(n)$.
\end{pf}

The analogous partition theoretic interpretation of~\eqref{rs1} is given next.
\begin{thm}\label{a1}
Let $A_1(n)$ denote the number of partitions of $n$ such that if
$2j$ is the largest repeated even part, then all positive even integers less than $2j$
also appear at least twice, no odd part less than $2j+2$ appears, and
no part greater than $2j$ is repeated.  Then $A_1(n) = C_{3,1}(n)$ for
all $n$.
\end{thm}

\begin{pf}
The proof parallels that of Theorem~\ref{andrews}, except that the role of 
$(-q^{2j+1};q)_\infty$ in Theorem~\ref{andrews} is played by
$(-q^{2j+2};q)_\infty$, which is the generating function for partitions
into distinct parts with each part at least $2j+2$, and that the right hand side
of~\eqref{rs1} is the generating function for $C_{3,1}(n)$.
\end{pf}

The purpose of this paper is provide an explicit bijection between the partitions
enumerated by Andrews' $A_2(n)$ and those of Gordon's $G_{3,2}(n)$.
As we shall see, the same map also provides a bijection between
the partitions enumerated by $G_{3,1}(n)$ and those enumerated by $A_1(n)$,
as well as a new partition theorem related to a special case of 
David Bressoud's even modulus analog of Gordon's theorem.
 
 \section{Definitions}\label{Defs}
 \subsection{Standard Definitions in Partition Theory}
 The definitions and symbols introduced in this subsection are all 
standard (cf.~\cite{igm}).
 
 A \emph{partition} $\pi$ of an integer $n$ is a nonincreasing sequence 
of nonnegative integers
\[ \pi = \{ \pi_1, \pi_2, \pi_3, \dots \} \]
such that $\sum_{i=1}^\infty \pi_i = n$.
Each nonzero term in $ \{ \pi_1, \pi_2, \pi_3, \dots \} $ is called a \emph{part}
of the partition $\pi$.   The number of parts in $\pi$ is called the 
\emph{length} $\ell(\pi)$ of $\pi$.  
Since the tail $\{ \pi_{\ell(\pi)+1}, \pi_{\ell(\pi)+2}, \dots \}$ of any partition $\pi$
must be $\{ 0,0,0,0,\dots \}$, it will be convenient to suppress the infinite
string of zeros when writing a specific partition.
The \emph{empty partition}, $\emptyset = \{ 0,0,0,0,\dots \} = \{ \}$, has length zero,
i.e. no parts.

 For two partitions, $\pi$ and $\lambda$, 
we may write 
$\pi \geqq \lambda$ if $\pi_i\geqq \lambda_i$ for all $i$.
The \emph{multiplicity} of the integer $j$ in $\pi$, denoted
$m_j(\pi)$, is the number of times $j$ appears in $\pi$.  
 
  At times it will be convenient to express $\pi$ in the alternate notation
  \[\pi = \langle 1^{m_1(\pi)} 2^{m_2(\pi)} 3^{m_3(\pi)} \cdots \rangle \]
meaning that $\pi$ contains $m_1(\pi)$ 1's, $m_2(\pi)$ 2's, etc.  In this notation it
is customary to omit the term $j^{m_j(\pi)}$ when $m_j(\pi)=0$.
  
 The \emph{union} of two partitions $\pi$ and $\lambda$, denoted $\pi\cup\lambda$,
 is the partition
 whose parts are those of $\pi$ and $\lambda$ together, arranged in nonincreasing order.
 For example,
   \[ \{ 8,3,3,2,1\} \cup \{ 9, 7, 5, 3, 1 \} = \{ 9,8,7,5,3,3,3,2,1,1 \} .\]
   
 The \emph{sum} of two partitions $\pi$ and $\lambda$ is
 \[ \pi + \lambda := \{ \pi_1+\lambda_1, \pi_2+\lambda_2, \dots \}. \]
 If $\pi\geqq\lambda$, then the \emph{difference} of $\pi$ and $\lambda$ is given by
 \[ \pi-\lambda :=\{ \pi_1-\lambda_1, \pi_2-\lambda_2, \dots \} .\]

 \subsection{Definitions of Special Symbols}
 In the standard notation, the $i$th largest part of the partition $\pi$ is denoted $\pi_i$, where
 $1\leqq i \leqq \ell(\pi)$.
 Accordingly, the $i$th smallest part of $\pi$ is $\pi_{\ell(\pi)-i+1.}$  However, there will
 be occasions where a less cumbersome notation for the $i$th smallest part will be
 useful, so let us define
 \[ \pi_{[i]} := \pi_{\ell(\pi)-i+1} \mbox{\ for $i=1,2,\dots,\ell(\pi).$} \]

  For two partitions $\pi$ and $\lambda$, let us write
 $\pi\succ\lambda$ if the smallest part of $\pi$ is greater than the largest part of $\lambda$, i.e.,
 \[ \pi_{[1]} = \pi_{\ell(\pi)}   > \lambda_1. \]
 Let $D(\pi)$ denote the number of different parts in
 $\pi$ which appear
at least twice, i.e.
\[ D(\pi):= \card\{  j \in \pi \ | \ m_j(\pi)>1 \}, \]   
thus $D(\{21, 15,15,12, 11, 9,9, 6, 5,5, 2 \}) = 3$ 
since three integers (15, 9, and 5) each appear
more than once as parts.
 
  Let $R_k(\pi)$ denote the $k$-th largest repeated part in $\pi$, i.e.
  \[ R_1(\pi):= \max\{ j \ | \ m_{j}(\pi) > 1\}, \]
  \[ R_k(\pi):= \max\{ j \ | \ m_{j}(\pi)>1 \mbox{ and  } j< R_{k-1} (\pi) \} \mbox{ for $k=2,3,\dots, D(\pi)$,}\]
  \[ R_{D(\pi)+1} (\pi) = 0. \]

  Let $\Gone$ denote the set of partitions enumerated by $G_{3,1}(n)$ in
Gordon's theorem, i.e. partitions $\pi$ such that
 \begin{equation} m_{1}(\pi) = 0 \label{g1cond1}\end{equation} and
 \begin{equation} m_{j} (\pi) + m_{j+1}(\pi) \leqq 2 \label{g1cond2}\end{equation} for all $j\geqq 1$.
 
  Let $\Gtwo$ denote the set of partitions enumerated by $G_{3,2}(n)$ in
Gordon's theorem, i.e. partitions $\pi$ such that
 \begin{equation} m_{1}(\pi) \leqq 1 \label{g2cond1}\end{equation} and
 \begin{equation} m_{j} (\pi) + m_{j+1}(\pi) \leqq 2 \label{g2cond2}\end{equation} for all $j\geqq 1$.
 
  Let $\Aone$ denote the set of partitions enumerated by $A_1(n)$ in
Theorem~\ref{a1}, i.e. partitions $\pi$ such that
\begin{equation} m_j(\pi) \leqq 1 \mbox{  if $j$ is odd,} \label{a1cond1} \end{equation}
\begin{equation} m_j(\pi)  = 0 
   \mbox { if $j$ is odd and $j<R_1(\pi)+2$, and}\label{a1cond2} \end{equation}
\begin{equation} m_j(\pi)\geqq  2  
    \mbox{ if $j$ is even and $j<R_1(\pi)$.}\label{a1cond3} \end{equation}

  Let $\Atwo$ denote the set of partitions enumerated by $A_2(n)$ in
Theorem~\ref{andrews}, i.e. partitions $\pi$ such that
\begin{equation} m_j(\pi) \leqq 1 
    \mbox{  if $j$ is odd,} \label{a2cond1}  \end{equation}
\begin{equation} m_j(\pi)  = 0 
   \mbox{ if $j$ is odd and $j<R_1(\pi)$, and}\label{a2cond2} \end{equation}
\begin{equation} m_j(\pi)\geqq  2  
    \mbox{ if $j$ is even and $j<R_1(\pi)$.}\label{a2cond3} \end{equation}

 For any partition $\pi\in\Gtwo$, let 
 $\pi^{(0)}$ denote the subpartition of $\pi$ whose parts are greater than $R_1(\pi)$,
 and $\pi^{(k)}$ denote the subpartition of $\pi$ whose parts are less than $R_k (\pi)$ and
 greater than $R_{k+1}(\pi)$, for $1\leqq k\leqq D(\pi)$.  
 Notice that by Condition~\eqref{g2cond2}, no part of $\pi\in\Gtwo$
 may appear more than twice.  Thus any $\pi\in\Gtwo$ may be decomposed 
uniquely as
 \[ \pi = \pi^{(0)} \cup  \langle R_1(\pi)^2 \rangle \cup  
   \pi^{(1)} \cup  \langle R_2(\pi)^2 \rangle
 \cup \pi^{(2)} \cup
 \cdots  \cup \langle R_{D(\pi)}(\pi)^2 \rangle \cup \pi^{(D(\pi))} ,\]
 where each $\pi^{(k)}$ is a partition into distinct parts, and
 \[ \pi^{(0)} \succ  \langle R_1(\pi)^2 \rangle \succ  
   \pi^{(1)} \succ  \langle R_2(\pi)^2 \rangle
 \succ \pi^{(2)} \succ
 \cdots  \succ \langle R_{D(\pi)}(\pi)^2 \rangle \succ 
  \pi^{(D(\pi))} .\]  Of course, any of the
 $\pi^{(k)}$ may be empty.
 
  Let $P(\pi)$ denote the subpartition of $\pi$ consisting of parts greater than $R_1(\pi)$, and
 $S(\pi)$ denote the subpartition of $\pi$ consisting of parts less than $R_1(\pi)$. 
 Thus, $$P(\{21, 15,15,12, 11, 9,9, 7, 5,5, 2 \}) = \{ 21 \}$$ and
 $$S(\{21, 15,15,12, 11, 9,9, 7, 5,5, 2 \}) = \{12, 11, 9,9, 7, 5,5, 2 \} .$$

  For any partition $\lambda\in\Atwo$, $\lambda$ can be decomposed uniquely
into 
 \[ \lambda = P(\lambda) \cup L(\lambda), \]
where  $L(\lambda)$
is the subpartition of $\lambda$ consisting of all parts of $\lambda$
less than or equal to $R_1(\lambda)$.  Clearly, $P(\lambda)$ is
a partition with distinct parts, $R_1(\lambda)$ is a partition
into parts which are even and repeated, and
\[ P(\lambda) \succ L(\lambda). \]

\begin{rem}
Note that $P(\lambda) = \lambda^{(0)}$. One reason for introducing 
the $P$ notation is to emphasize that for $\lambda\in\Atwo$,
the partition $\lambda$ naturally decomposes into just two subpartitions, 
$P(\lambda)$ and
$L(\lambda)$, whereas for $\pi\in\Gtwo$, the partition $\pi$
naturally decomposes into $2D(\pi)+1$ subpartitions, of which
all of the $\pi^{(k)}$ are partitions with distinct parts.  
Also, when $P$ is applied to a partition $\pi\in\Gtwo$,
it is used in conjunction with $S$, where $P$ and $S$ can be 
thought of as a ``prefix" and ``suffix" respectively.
\end{rem}

\section{The Bijection}
\subsection{Mapping $\Gone$ onto $\Aone$.}
Since $\Gone\subsetneqq\Gtwo$ and $\Aone\subsetneqq\Atwo$, 
let us initially turn our attention to $\Gone$ and $\Aone$.

\begin{defn} \label{fdef}Define the map $f$ on the set $\Gone$ by
\begin{equation} \label{fdefeq}
f(\pi) := \bigcup_{i=1}^{D(\pi)} \langle (2i)^{R_i(\pi) - R_{i+1}(\pi) - \ell(\pi^{(i)})}\rangle 
\cup \bigcup_{i=0}^{D(\pi)} \left( \pi^{(i)} + \langle  (2i)^{\ell(\pi^{(i)})}\rangle \right) ,\end{equation}
for $\pi\in\Gone$.
Equivalently, $f$ may be defined recursively by
  \[ f(\pi) := P(\pi) \cup \left[ \langle 2^{R_1(\pi)} \rangle + f( S(\pi) ) \right] \] with initial condition
  \[ f(\emptyset) = \emptyset. \]
\end{defn}

\begin{exmp}\label{ex1}
If
\[ \pi = 
\{40,37, 36, 22, 22, 20, 19, 17, 17, 15, 13, 12, 10, 8, 8, 4, 4, 2 \}, \]
then
   \begin{equation}
   f(\pi) = \{ 40, 37, 36, 22, 21, 19, 17, 16, 14, 10 \} \cup \langle 2^3 4^5 6^4 8^3\rangle . 
    \end{equation}
   
   In more detail, the mapping is
\begin{align*}& \qquad\{40,37, 36, 22, 22, 20, 19, 17, 17, 15, 13, 12, 10, 8, 8, 4, 4, 2 \} \\
&= \quad\{ 40 , 37, 36 \} \cup \{ 22,22 \} \cup \{ 20, 19 \} \cup  \{ 17,17 \} \cup \{ 15,13,12,10 \}\cup
\{8,8\} \cup \{ 4,4 \} \cup \{2 \}\\
&\overset{f}{\longrightarrow}  
\langle 2^{22-17-2} \rangle \cup \langle 4^{17-8-4}  \rangle\cup \langle 6^{8-4-0}\rangle \cup
\langle 8^{4-0-1} \rangle \\
&\qquad\qquad\cup
\{ 40 , 37, 36 \} \cup \{ 20+2, 19+2 \} \cup   \{ 15+4,13+4,12+4,10+4 \} \cup \emptyset
\cup \{2+8 \} \\
&=\quad \{ 40, 37, 36, 22, 21, 19, 17, 16, 14, 10 \} \cup \langle 2^3 4^5 6^4 8^3\rangle . 
\end{align*}

 The following diagram helps to demonstrate how the the recursive formulation of $f$ works.  
In each iteration, the largest repeated part $R$ is underlined, and then converted to
$R$ $2$'s in the next iteration.  At the last step, the columns are added to form the
image of $\pi$ under $f$.
 \[ \begin{array}{rrrrrrrrrr rrrrrrrr rrrrrrr}
40 & 37 & 36 & \underline{22} & \underline{22} & 20 & 19 & 17 & 17 & 15 & 13&12&10&8&8&4&4&2&&&&&&&
\end{array}
\] 
\[ \downarrow \]

  \[ \begin{array}{rrrrrrrrrr rrrrrrrrrrrrrrr}
40 & 37 & 36 & 2 & 2 & 2 & 2 & 2 & 2 & 2 & 2&2&2&2&2&2&2&2&2&2&2&2&2&2&2\\
      &       &    &20&19& \underline{17} & \underline{17} & 15 & 13 & 12 & 10&8&8&4&4&2
\end{array}
\]
\[ \downarrow \]
\[ \begin{array}{rrrrrrrrrr rrrrrrrrrrrrrrr}
40 & 37 & 36 & 2 & 2 & 2 & 2 & 2 & 2 & 2 & 2&2&2&2&2&2&2&2&2&2&2&2&2&2&2\\
      &       &     &20&19& 2 & 2 & 2 & 2 & 2 & 2&2&2&2&2&2&2&2&2&2&2&2&  & &\\
      &        &     &      &    &15&13&12&10& \underline{8} & \underline{8}&4&4&2& & & 
     \end{array}
\] \[ \downarrow \]
\[ \begin{array}{rrrrrrrrrr rrrrrrrrrrrrrrr}
40 & 37 & 36 & 2 & 2 & 2 & 2 & 2 & 2 & 2 & 2&2&2&2&2&2&2&2&2&2&2&2&2&2&2\\
      &       &     &20&19& 2 & 2 & 2 & 2 & 2 & 2&2&2&2&2&2&2&2&2&2&2&2&  & &\\
      &        &     &      &    &15&13&12&10& 2 & 2&2&2&2&2&2&2\\
      &        &      &     &      &    &     &     &   & \underline{4} &\underline{ 4} & 2 &  & \\
\end{array}
\]
\[ \downarrow \]
\begin{equation}\begin{array}{rrrrrrrrrr rrrrrrrrrrrrrrr}
40 & 37 & 36 & 2 & 2 & 2 & 2 & 2 & 2 & 2 & 2&2&2&2&2&2&2&2&2&2&2&2&2&2&2\\
      &       &  +  &20&19& 2 & 2 & 2 & 2 & 2 & 2&2&2&2&2&2&2&2&2&2&2&2&  & &\\
      &        &     &      &   +&15&13&12&10& 2 & 2&2&2&2&2&2&2\\
      &        &      &     &      &    &     &     &  +& 2 & 2 & 2 & 2\\
       &        &     &      &     &   &  &  & + & 2 \\
       \hline
40 & 37 & 36& 22& 21& 19 & 17 & 16 & 14& 10 & 8 &8&8&6&6&6&6&4&4&4&4&4&2&2&2
\end{array}
\end{equation}
\end{exmp}

Thus the parts of both $\pi$ and $f(\pi)$ are encoded
in the matrix
\begin{equation}\left[ \begin{array}{rrrrrrrrrr rrrrrrrrrrrrrrr}
40 & 37 & 36 & 2 & 2 & 2 & 2 & 2 & 2 & 2 & 2&2&2&2&2&2&2&2&2&2&2&2&2&2&2\\
   0   &  0     &  0  &20&19& 2 & 2 & 2 & 2 & 2 & 2&2&2&2&2&2&2&2&2&2&2&2&0  &0 &0\\
    0  &   0     &   0  &  0    &   0&15&13&12&10& 2 & 2&2&2&2&2&2&2&0&0&0&0&0&0&0&0\\
    0  &    0    &   0  &   0  &  0    &    0&  0   & 0    &  0& 2 & 2 & 2 & 2&0&0&0&0&0&0&0&0&0&0&0&0\\
      0 &   0     &    0 & 0    &  0   & 0  & 0 & 0 & 0 & 2 &0&0&0&0&0&0&0&0&0&0&0&0&0&0&0
\end{array} \right] . \label{fdiag-ex}
\end{equation}

Display~\eqref{fdiag-ex} may be the best way of simultaneously displaying $\pi$ and $f(\pi)$;
the distinct parts of $\pi$ appear in the matrix as nonzero, nonrepeated entries in a
given row, parts
$R_i(\pi)$ which
appear twice are represented as a sequence of $R_i(\pi)$ $2$'s in the $i$th row, and
the parts of $f(\pi)$ are the sums of the columns.

 Accordingly,~\eqref{fdiag-ex} motivates the following definition.  Let $\lambda:= f(\pi)$.
 \begin{defn} Let $\pi\in\Gtwo$.  
The \emph{$\mathcal{S}$-diagram of the partitions $\pi$ and $\lambda$} 
is the $ \left( D(\pi)+1 \right) \times \left( \ell(\pi^{(0)}) + R_1(\pi) \right)$ matrix 
(or, equivalently, the $\left( 1+R(\lambda)/2 \right) \times \ell(\lambda)$ matrix) whose
first row consists of the parts of $\pi^{(0)}$ in nonincreasing order 
followed by $R_1(\pi)$ $2$'s, and whose 
$i$th row consists of $\sum_{k=0}^{i-1} \ell( \pi^{(k)} ) \ 0$'s followed by
the parts of $\pi^{(i-1)}$ in nonincreasing order followed by $R_i(\pi)$ $2$'s, followed by
$0$'s, for $2\leqq i \leqq D(\pi)+1 $.  
The parts of $\lambda$ are then given by the sums of the columns.
\end{defn}

Clearly, there is a unique $\mathcal{S}$-diagram for each partition $\pi\in\Gtwo$.
Next, let us examine how to determine the $\mathcal{S}$-diagram given only a
partition in $\Atwo$.

 Again, $\lambda=f(\pi)$.  Observe that 
\begin{equation*}
  \lambda_j = \pi_j, \mbox{\  if $1\leqq j \leqq \ell(\pi^{(0)})$ },
\end{equation*}
\[ \ell(\lambda) = \ell( \pi^{(0)} ) + R_1(\pi) ,\]
and \[ \lambda_1 > \lambda_2 > \dots > \lambda_{\ell(\pi^{(0)})} \] since
$\pi^{(0)}$ is a partition with distinct parts.
Accordingly, 
\[ \ell(\lambda) - j + 1 < \lambda_j, \mbox{ \ if $1\leqq j \leqq \ell ( \pi^{(0)}) $} \]
However, \[  \lambda_{\ell(\pi^{(0)} ) } = 2 + \pi_{\ell(\pi^{(0)}) + 3} \leqq 2+2+ R_1(\pi) \]
and so \[ \ell{(\lambda)} - (\ell(\pi^{(0)} + 1)) + 1 - \lambda_{\ell(\pi^{(0)} + 1)} \geqq   \ell(\pi^{(0)} + 1) \]
which means that
\begin{equation}\label{R1}
  R_1(\pi) = \ell{(\lambda)} - (\ell(\pi^{(0)} + 1)) + 1 - \lambda_{\ell(\pi^{(0)} + 1)}
  = \max_{1\leqq j \leqq \ell(\lambda)}
  \left\{  j \ \Big| \ j \geqq \lambda_{[j]} \right\}.
\end{equation}
In order to clarify what is being asserted, let us examine the above statements
in terms of Example~\ref{ex1}.  The three largest parts of $\lambda=f(\pi)$, namely 40, 37, and 36, 
are respectively the $25$th, $24$th, and $23$th \emph{smallest} parts.  These parts (40, 37, and 36),
are all larger than the number of parts less than equal to them ($25$, $24$, and $23$ respectively).
However the fourth largest (i.e. $22$nd smallest) part of $\lambda$, $\lambda_4=\lambda_{[22]}=22$ is such that
it has $22$ parts less than or equal to it, and since the part in question, $22$, does
not exceed the number of parts less than or equal to it (i.e. $22$), it must be
that $R_1(\pi) = 22$.
  
   In other words, knowing the parts of $\lambda$ (with no foreknowledge of the
 parts of $\pi$),~\eqref{R1} allows us to 
 recover the first row of the $\mathcal{S}$-diagram.

   Analogously, knowing the parts of $\lambda$ only, we can recover the $i$th row of the $\mathcal{S}$-diagram
via
\begin{equation}\label{R}  
 R_i(\pi) = \max_{1\leqq j \leqq \ell(\lambda)}
  \left\{ j - \sum_{h=1}^{i-1} m_{2h}(\lambda)\ \Big| \ 
  j - \sum_{h=1}^{i-1} m_{2h}(\lambda)\geqq \lambda_{[j]}-2(i-1) 
  \right\}
\end{equation}
for $i = 1,2,\dots, 1+\frac{R_1(\lambda)}{2}.$
   
\begin{thm}\label{G1toA1}
The map $f$ is a bijection of $\Gone$ onto $\Aone.$
\end{thm}

\begin{pf}
First, we need to show that for any $\pi\in\Gone$, $f(\pi)\in\Aone$.
Let $1\leqq k \leqq D(\pi)$.
Since $R_{k+1} (\pi)$ appears twice in $\pi$, the smallest part in $\pi^{(k)}$
must be at least two less than $R_{k+1}(\pi)$ (by~\eqref{g2cond2}), and
thus the largest part in $\pi^{(k+1)}$ must be at least two larger than $R_{k+1}(\pi)$,
(again by~\eqref{g2cond2}).  Thus, the smallest part in $\pi^{(k)}$ must
exceed the largest part in $\pi^{(k+1)}$ by at least four.
Therefore,
\[ \left( \pi^{(k)} + \langle (2k)^{\ell(\pi^{(k)})} \rangle \right) \succ 
\left( \pi^{(k+1)} + \langle (2k+2)^{\ell(\pi^{(k+1)})} \rangle \right), \]
and thus all parts of $P(f(\pi))$ are distinct.  Furthermore,
\[ P(f(\pi)) \succ L(f(\pi)) \]
since
\begin{align*} P(f(\pi))_{[1]} &= \pi_{[1]} +2D(\pi)\\
                                                                &=\pi_{[1]} + L(f(\pi))_1\\
                                                                 &\geqq 2 + L(f(\pi))_1\qquad \mbox{  (by~\eqref{g2cond1})}\\
                                                                 &> L(f(\pi))_1.
\end{align*} Also, by~\eqref{fdefeq}, $L(f(\pi))$ contains the numbers $2i$ for $i=1,2,\dots,D(\pi)$
as parts. 
The number of appearances of $2i$ in $f(\pi)$ is $R_i(\pi) - R_{i+1}(\pi) - \ell( \pi^{(i)})$.
  Notice that Condition~\eqref{g1cond2} forces
$R_i(\pi) - R_{i+1}(\pi) - \ell( \pi^{(i)}) \geqq 2$, and so $f(\pi)\in\Aone$.

  Next, we demonstrate the invertibility of $f$.
Again, let $\lambda=f(\pi)$.
Notice that $L(\lambda)_1 = R_1(\lambda) = 2 D(\pi)$
and that $\ell(P(\lambda)) = \sum_{k=0}^{D(\pi)} \ell(\pi^{(k)})$.
From the definition~\eqref{fdef} of $f(\pi)$,
  \begin{equation}\label{Rl}
 m_{2j}(\lambda) = R_j(\pi)- R_{j+1}(\pi) - \ell( \pi^{(j)} )
   \end{equation}
for $j=1,2,\dots,D(\pi)$. 
  Thus, once the $R_j(\pi)$ and $\ell(\pi^{(j)})$ are recovered
from $\lambda$ using~\eqref{R}, 
this will be sufficient information to reconstruct
$\pi = f^{-1}(\lambda)$ from~\eqref{fdefeq}.

Once the $R_j(\pi)$ are known, the rest of $\pi$ can be easily 
constructed by subtracting appropriate multiples of two from
the parts of $P(\lambda)$ with the aid of~\eqref{fdefeq} and~\eqref{Rl}. 
\end{pf}

\begin{exmp}
Let us now consider $\lambda=
\{40,37,36,22,21,19,17,16,14,10\}\cup \langle 2^3 4^5 6^4 8^3\rangle\in\Aone$ 
in order to show that it will map to the $\pi$ from Example~\ref{ex1}.

  Since the initial repeated even parts
in $\lambda$ are $2,4,6,$ and $8$, a total of 4 different parts,
$\pi$ must contain $4$ different repeated parts,
$R_1(\pi)> R_2(\pi)>R_3(\pi) > R_4(\pi)$.  
  Other important quantities which can be simply ``read off"
are $\ell(\lambda)=25$, $m_2(\lambda)=3$, 
$m_4(\lambda)=5$, $m_6(\lambda)=4$, and $m_8(\lambda)=3$.
Now consider the computation of $R_1(\pi)$ using~\eqref{R}. 
The $25$ parts of $\lambda$ are labeled from smallest to largest.
\[
\begin{array}{rrrrrrrrrrrrrrrrrrrrrrrrr}
 25 & 24 &23 &22 &21 &20 & 19 & 18 &17 &16 &15 &14 &13 &12 &11 &10 &9 &8 &7
    &6 & 5& 4&3 & 2 &1\\
 40 &37 & 36 & 22 & 21 & 19 & 17 & 16 & 14 & 10 & 8&8&8&6&6&6&6&4&4&4&4&4&2&2&2
\\
 & & &\uparrow&& &&&&&&&&&&&&&&&&& & &
\end{array}
\] Reading from the left, the first instance of a part not exceeding the number of parts
less than or equal to it is $\lambda_{[22]} = 22$,
  thus, $R_1(\pi) = 22$, and $\pi^{(0)} = \{ 40,37,36 \}$.
Therefore, the first row of the $\mathcal S$-diagram must be
$40$, $37$, $36$, followed by twenty-two $2$'s.

Since $\lambda$ has $25$ parts, three of which are $2$'s, the partition
$\lambda - \langle 2^{25} \rangle $ contains $25-3 = 22$ parts.  
Furthermore, the three largest parts of $\lambda$ have already been
determined in the previous step, so they can be removed from
further consideration.  Labeling 
the $19$ parts of $\{ \lambda_4, \lambda_5, \dots, \lambda_{22} \} - \langle 2^{19} \rangle $  from smallest to largest,
\[
\begin{array}{rrrrrrrrrrrrrrrrrrr}
 19 & 18 &17 &16 &15 &14 &13 &12 &11 &10 &9 &8 &7
    &6 & 5& 4&3 & 2 &1\\
 20 & 19 & 17 & 15 & 14 & 12 & 8 & 6&6&6&4&4&4&4&2&2&2&2&2 
\\
 & &\uparrow &&&&&&&&&&&&&& & & 
\end{array}
\]
which lets us conclude that $R_2(\pi) = 17$ and $\pi^{(1)} = \{ 22-2, 21-2 \}$. Therefore, we now know that the
second row of the $\mathcal S$-diagram begins with $\ell(\pi^{(0)}) = 3$ zeros, followed by
$22-2, 21-2$, followed by seventeen $2$'s, followed by $m_2(\lambda)=3$ zeros.

To find the next row, number the parts of $ \{ \lambda_6, \lambda_7, \dots, \lambda_{17} \} - 
\langle 4^{12} \rangle $
 from smallest to largest:
\[
\begin{array}{rrrrrrrrrrrrr}
  12 & 11 &10  &   9 & 8 &7 & 6 & 5 & 4 &3 &2&1 \\
 15 & 13 & 12 & 10 & 6 & 4& 4 & 4 & 2 &2 &2&2\\
              &       &      & &\uparrow &    &    &    &     &   &   &
\end{array}
\]
Thus, $R_3(\pi) = 8$ and $\pi^{(2)} = \{ 19-4, 17-4, 16-4, 14-4 \}$. 

Next,  
\[
\begin{array}{rrrr}
   4 &3 &2&1 \\
  4 & 2& 2 & 2 \\
 \uparrow    &    &    &   
\end{array}
\]
and so $R_4(\pi)=4$ and $\pi^{(3)} = \emptyset$.

Finally, all that remains is the single part $10- 8 = 2$, so this $2$ must be placed in column
$9 =1+ \ell(\pi^{(0)}) + \ell(\pi^{(1)}) + \ell(\pi^{(2)}) + \ell(\pi^{(3)})$ of the last row
of the $\mathcal S$-diagram.

Thus, the full $\mathcal S$-diagram is given by~\eqref{fdiag-ex}, and we conclude that
 \[ \pi = f^{-1}(\lambda) = \{ 40,37,36, 22,22, 20,19, 17,17, 15, 13,12,10, 8,8, 4,4, 2 \}. \]
\end{exmp}

\subsection{Extending the Map to $\Gtwo\to\Atwo$}
\begin{thm}
The map $f$ given in Definition~\ref{fdef} also provides a bijection from
$\Gtwo$ to $\Atwo$.
\end{thm}
\begin{pf}
Notice that $\Gone\subsetneqq\Gtwo$ and $\Aone\subsetneqq\Atwo$.
The set $\Gtwo\setminus\Gone$ consists of precisely 
those partitions in $\Gtwo$ which
contain a one. The set $\Atwo\setminus\Aone$ consists of precisely those
partitions $\lambda$ in $\Atwo$ for which the smallest part
of $P(\lambda)$ is exactly one more than the largest part in 
$L(\lambda)$.   Observe that any partition $\pi\in\Gtwo\setminus\Gone$ has
an $\mathcal{S}$-diagram whose 
$\left( R(\pi)+1, \sum_{h=0}^{D(\pi)} \ell(\pi^{(h)}) \right)$ entry is a $1$.
This condition is equivalent to the sum of the 
$\left( \sum_{h=0}^{D(\pi)} \ell(\pi^{(h)}) \right)$th
column being exactly one more than the next column, i.e. the smallest part of
$P(\lambda)$ is exactly one more than the largest part of $L(\lambda)$.
\end{pf} 

\subsection{The Other Rogers-Selberg Identity~\eqref{rs3} }
Obviously, it would be desirable to use the $f$ map, or some generalization of it, to
produce a bijection from the set of partitions $\Gthree$ enumerated 
by the $G_{3,3}(n)$ of Gordon's theorem
onto a set of partitions enumerated by the left hand side of~\eqref{rs3}, say $\Athree$.  
In particular, we would like to have $\Atwo\subsetneqq\Athree$, just as 
we have $\Aone\subsetneqq\Atwo$.

Difficulties arise immediately, however.  
Firstly, the unaltered $f$-map will not
produce a bijection from the partitions enumerated by $G_{3,3}(n)$; observe that
$f(\{3,1,1\}) = f(\{3,2\}) = \{3,2\}$.  
Furthermore, finding a natural partition theoretic interpretation
of the left hand side of~\eqref{rs3} analogous to that of Theorem~\ref{andrews} is
not as straightforward as one might first suppose.  
In~\eqref{rs1} and~\eqref{rs2}, the 
$2j^2+2j = 2+2+4+4+\cdots+2j+2j$ in the exponent of $q$ over the 
expression $(q^2;q^2)_j$ allows
the expression \[ \frac{q^{2j^2+2j} }{(q^2;q^2)_j } \] 
to be dealt with neatly as partitions into
$2$'s, $4$'s, \dots, $2j$'s with each part appearing at least twice.  
In~\eqref{rs3}, we
instead have $q^{2j^2},$ but how shall the $2j^2$ be split up?  

 We will allow ourselves to be guided by a simple bijection of $\Gthree$ onto $\Gone$.
  Consider the map 
   \[ g: \Gthree \to \Gone \] 
where $g(\{ \pi_1, \pi_2, \dots , \pi_{\ell(\pi)} \}) = \{ \pi_1 + 1, \pi_2 + 1, \dots, \pi_{\ell(\pi)}+1 \} $.
It is clear that $\ell(\pi) = \ell(g(\pi))$ and that
if $\pi$ is a partition of $n$, then $g(\pi)$ is a partition of $n+\ell(\pi)$.  It is also clear
that $g$ is a bijection.  By Theorem~\ref{G1toA1}, $f$ maps $\Gone$ bijectively onto
$\Aone$. Thus, all that is needed is to map $\Aone$  bijectively onto $\Athree$ in as
simple a manner as possible, so that a partition $f(g(\pi))$ of $n+\ell(\pi)$ in  $\Aone$
maps to a partition of $n$ in $\Athree$.   
  Define $h : \Aone \to \Athree$ as follows. 
For $\lambda\in\Aone$, subtract 1 from each part in $P(\lambda )$ and subtract 1 from two 2's, two 4's, two 6's, . . . , two ($R_1(\lambda)$)'s.
Letting $\bar{f}:= h\circ f\circ g$, we have a bijection
\[ \bar{f} : \Gthree \to \Athree .\] 
To view $\bar{f}$ as a single operation, rather than a composition of three maps, let us
define the following modified $\mathcal{S}$-diagram.

  \begin{defn} Let $\pi\in\Gthree$ and $\lambda=\bar{f}(\pi)$. 
The \emph{$\bar{\mathcal{S}}$-diagram of the partitions $\pi$ and $\lambda$} 
is the $ \left( D(\pi)+1 \right) \times \left( \ell(\pi^{(0)}) + R_1(\pi) + 1 \right)$ matrix 
(or, equivalently, the $\left( 1+R(\lambda)/2 \right) \times \ell(\lambda)$ matrix) whose
first row consists of the parts of $\pi^{(0)}$ in nonincreasing order 
followed by $R_1(\pi)-1$ $2$'s, then two $1$'s and whose 
$i$th row consists of $\sum_{k=0}^{i-1} \ell( \pi^{(k)} ) \ 0$'s, followed by
the parts of $\pi^{(i-1)}$ in nonincreasing order, followed by $R_i(\pi)-1$ $2$'s, followed by
two $1$'s, and the rest
$0$'s, for $2\leqq i \leqq D(\pi)+1 $.  
The parts of $\lambda$ are then given by the sums of the columns.
\end{defn}

\begin{exmp} $\bar{f} (\{ 16,14,12,12,7,5,5,3,2,1 \} ) =
\{ 16,14,9,7,6,5,4,3,3, 2,2,2,2,1,1 \} .$
The corresponding $\bar{\mathcal S}$-diagram is
\[ \left[  
  \begin{array}{rrrrrrrrrrrrrrr}
  16 & 14 & 2& 2& 2& 2&  2 & 2 & 2 & 2 & 2 & 2 & 2& 1 & 1 \\
   0   &  0  & 7 & 2 & 2& 2 & 2  & 1 & 1 & 0 & 0 & 0 & 0 & 0 & 0 \\
   0  &  0 & 0 & 3 & 2 & 1 & 0 & 0 & 0 & 0 & 0 & 0 & 0 & 0 & 0 
   \end{array} \right].
  \]
\end{exmp}

The partition theorem analogous to Theorems~\ref{andrews} and \ref{a1} is therefore
as follows:
\begin{thm}
Let $A_3(n)$ denote the number of partitions of $n$ such that if $2j-1$ is the largest
repeated odd part, then all positive odd integers less than $2j$ appear exactly twice,
and no part greater than $2j$ may be repeated.  Then $A_3(n) = C_{3,3}(n)$ for
all $n$.
\end{thm}
\begin{pf}
Note that
\[ \frac{q^{2j^2}(-q^{2j+1};q)_\infty}{(q^2;q^2)} 
 = \frac{ q^{1+1+3+3+\cdots+(2j-1)+(2j-1)} }{(q^2;q^2)} \times (-q^{2j+1};q)_\infty.
\]
By the methods of Euler (cf.~\cite[p. 4 ff.]{top}), the expression 
\[ \frac{ q^{1+1+3+3+\cdots+(2j-1)+(2j-1)} }{(q^2;q^2)} \]
is the generating function for partitions into exactly two $1$'s, two $3$'s, two $5$'s, $\dots$, two $(2j-1)$'s with $2$'s, $4$'s, \dots, $2j$'s allowed to appear any number of times (or not at all).
The product $(-q^{2j+1};q)_\infty$ is the generating
function for partitions into distinct parts with each part at least $2j+1$.  Thus, by summing over
all nonnegative $j$, it follows that the left hand side of~\eqref{rs3} is the
generating function for $A_3(n)$.  Again by Euler's method, it is 
immediate that
the right hand side of~\eqref{rs3} is the generating function for $C_{3,3}(n)$.
\end{pf}

\begin{rem}
 It should be noted that Andrews gave a different partition-theoretic interpretation of
the Rogers-Selberg identities
in~\cite{gea:gfp}.  This interpretation was studied further by Bressoud~\cite{dmb}.
\end{rem}

\section{The Modulus 6 Case of Bressoud's Theorem}
Bressoud~\cite{dmb:even} gave an analog of Gordon's partition theorem for even moduli.
\begin{thm}[Bressoud] \label{bressoud}
Let $B_{k,i}(n)$ denote the number of partitions $\pi$ of $n$ into parts such that 
  \begin{align}
      & m_1(\pi) \leqq i-1, \label{Bcond1}\\
      & m_j(\pi) + m_{j+1}(\pi) \leqq k-1 \mbox{ for } j=1,2,3,\dots, \mbox{ and } \label{Bcond2}\\
     & \mbox{if  } \pi_j - \pi_{j+k-2} \leqq 1, \mbox{ then  }\sum_{h=0}^{k-2} \pi_{j+h} \equiv (i-1)\pmod{2}.
     \label{Bcond3}
  \end{align}
    Let $D_{k,i}(n)$ denote the number of partitions of $n$
into parts $\not\equiv 0,\pm i \pmod{2k}$.  Then $B_{k,i}(n) = D_{k,i}(n)$ for 
$1\leqq i < k$ and all integers $n$.
\end{thm}
  The $k=3,\ i=1$ case of Bressoud's theorem first appeared in Andrews~\cite[p. 432, Thm. 1]{gea:snpt}.
This will be the case of most interest to us.   
It may be stated as follows.
\begin{thm}[Andrews]
The number of partitions of $n$ into parts greater than $1$, in which no consecutive integers
appear,
and in which no part appears more than twice equals the number of partitions of $n$ into
parts $\not\equiv 0, \pm 1\pmod{6}$. 
\end{thm}

The ``mod 6 analog'' of the Rogers-Selberg identity~\eqref{rs1} may be stated as
\begin{equation} \label{sl27}
 \sum_{j=0}^\infty \frac{ q^{2j^2+2j} (q;q^2)_j (-q^{2j+2};q)_\infty }{(q^2;q^2)_j}
 = \underset{j\not\equiv 0,\pm 1{\hskip-3mm}\pmod 6}{\prod_{j\geqq 1}}
 \frac{1}{1-q^j}.
\end{equation} 
Equation~\eqref{sl27} is due to Slater~\cite[p. 154, Eq. (27), with $q$ replaced by $-q$]{ljs}.

 The $k=3,\ i=2$ case of Bressoud's theorem
corresponds to Euler's classic theorem that the number of partitions into odd parts
equals the number of partitions into distinct parts~\cite[p. 5, Cor. 1.2]{top}, and so
the ``mod 6 analog" of the Rogers-Selberg identity~\eqref{rs2} is
simply

\begin{equation} \label{sl27-b}
 (-q;q)_\infty 
 = \underset{j\not\equiv 0,\pm 2{\hskip-3mm}\pmod 6}{\prod_{j\geqq 1}}
 \frac{1}{1-q^j} \left( = \underset{j\equiv 1{\hskip-3mm}\pmod 2}{ \prod_{j=1}^\infty} \frac{1}{1-q^j} \right).
\end{equation}
Since the $f$ map restricted to partitions with distinct parts is the identity map, this case will
not be of interest here.

A third partner also exists, but the infinite product on the right hand side is not as neat as that of~\eqref{sl27} or~\eqref{sl27-b} because $-3\equiv 3\pmod{6}$:
\begin{equation} \label{sl27-c}
 \sum_{j=0}^\infty \frac{ q^{2j^2} (q;q^2)_j (-q^{2j+1};q)_\infty }{(q^2;q^2)_j}
 = \prod_{j\geqq 1} \frac{(1-q^{3j})(1-q^{6j-3})}{1-q^j}.
\end{equation}
An identity equivalent to~\eqref{sl27-c} was found independently by 
James McLaughlin
in a computer search~\cite{jm}.

Notice that for a given $k$ and $i$, the conditions restricting the appearance of parts
in partitions enumerated by $G_{k,i}(n)$ in Gordon's theorem are the same as those
restricting the appearance of parts in partitions enumerated by $B_{k,i}(n)$ in Bressoud's
theorem, except that Bressoud's theorem contains condition~\eqref{Bcond3}, while
Gordon's theorem does not.  Thus if $\Bone$ denotes the set of partitions enumerated
by $B_{3,1}(n)$ in Bressoud's theorem, it is immediate that 
  \[ \Bone \subsetneqq \Gone. \]
Therefore $f$ maps $B_{3,1}(n)$ onto some proper subset of $\Aone$.  It fairly
straightforward to see that the subset of $\Aone$ in question is those partitions
in which the distinct parts greater than the largest repeated even part are all
nonconsecutive integers.  Formally, the $f$ map supplies the following theorem.
\begin{thm}\label{Tthm}
Let $B_{3,1}(n)$ be as in Bressoud's theorem.  Let $T(n)$ denote
the number of partitions of $n$ into parts such that if $2j$ is the largest repeated
even part, then all positive even integers less than $2j$ also appear at least twice,
no odd part less than $2j+2$ appears, no part greater than $2j$ is repeated, and
no two consecutive integers appear.  Then $B_{3,1}(n) = T(n)$ for all $n$.
\end{thm}

While it must be admitted that the preceding theorem is not the most elegant identity in the
theory of partitions, it would nonetheless be desirable, in the present context, to be able to 
identify the partitions enumerated by $T(n)$ with the left hand side of~\eqref{sl27}.  Unfortunately,
the presence of the factor $(q;q^2)_j$ in the numerator of the left hand side of~\eqref{sl27}
(which is absent from the Rogers-Selberg identities) creates an inclusion-exclusion
situation; i.e. unlike Eqs.~\eqref{rs1},~\eqref{rs2}, and~\eqref{rs3}, Eq.~\eqref{sl27}
is \emph{not} positive term-by-term.  Accordingly, any hope of getting an 
immediate
result like that of Theorem~\ref{andrews} is dashed.
The author tried to find a direct combinatorial proof of Theorem~\ref{Tthm},
but was unable to do so.

Undeterred, we proceed by defining the following two-variable analogs of Eqs.~\eqref{sl27}--\eqref{sl27-c}.
  \begin{align}
    F_1 (a)&:= F_1(a;q):= \sum_{n=0}^\infty \frac{ a^{2n} q^{2n(n+1)} (q;q^2)_n (-aq^{2n+2};q)_\infty }{(q^2;q^2)_n} \label{F1}\\
    F_2(a) &:= F_2(a;q) := (-aq;q)_\infty  \label{F2} \\
    F_3 (a)&:= F_3(a;q):= \sum_{n=0}^\infty \frac{ a^{2n} q^{2n^2} (q;q^2)_n (-aq^{2n+1};q)_\infty }{(q^2;q^2)_n}. \label{F3}
  \end{align}
Next, define
 \begin{align}
  E_1 (a)&:= E_1(a;q):= \sum_{n=0}^\infty \frac{ a^{2n} q^{2n(n+1)}  }{(q^2;q^2)_n}
   \sum_{m=0}^\infty \frac{ a^m q^{m(m+2n+1)}} {(q;q)_m} \label{E1}\\
    E_2(a) &:= E_2(a;q) := (-aq;q)_\infty \label{E2} \\
     E_3 (a)&:= E_3(a;q):= \sum_{n=0}^\infty \frac{ a^{2n} q^{2n^2}  }{(q^2;q^2)_n}
   \sum_{m=0}^\infty \frac{ a^m q^{m(m+2n)}} {(q;q)_m}. \label{E3}
  \end{align}
  
\begin{lem} The functions $F_1(a)$, $F_2(a)$, $F_3(a)$ satisfy the following
system of $q$-difference equations:
\begin{align}
  F_1(a) &= F_3(aq) \label{Feq1}\\
  F_2(a) &= (1+aq) F_2(aq) \label{Feq2} \\
  F_3(a) &= F_1(a) + aq(1+aq) F_1(aq) \label{Feq3}
\end{align}
\end{lem}
\begin{pf}
 Equations~\eqref{Feq1} and~\eqref{Feq2} are immediate from the 
definitions~\eqref{F1}--\eqref{F3}.
\begin{align*}
F_3(a) - F_1(a) &= \sum_{n=0}^\infty 
   \frac{a^{2n} q^{2n^2} (q;q^2)_n (-aq^{2n+2};q)_\infty}{(q^2;q^2)_n}
  (1+aq^{2n+1} - q^{2n} )\\
  &=\sum_{n=0}^\infty 
   \frac{a^{2n} q^{2n^2} (q;q^2)_n (-aq^{2n+2};q)_\infty}{(q^2;q^2)_n}
 \\ & \qquad\qquad\times
  \Big( (1-q^{2n})(1+aq^{2n+1}) + aq^{4n+1} \Big)\\
  &=\sum_{n=1}^\infty 
   \frac{a^{2n} q^{2n^2} (q;q^2)_n (-aq^{2n+1};q)_\infty}{(q^2;q^2)_{n-1} }
   \\ &\qquad\qquad +\sum_{n=0}^\infty 
   \frac{a^{2n+1} q^{2n^2+4n+1} (q;q^2)_n (-aq^{2n+2};q)_\infty}{(q^2;q^2)_n}\\
  &=\sum_{n=0}^\infty 
   \frac{a^{2n+2} q^{2n^2+4n+2} (q;q^2)_{n+1} (-aq^{2n+3};q)_\infty}{(q^2;q^2)_{n} }\\ &\qquad\qquad
   +\sum_{n=0}^\infty 
   \frac{a^{2n+1} q^{2n^2+4n+1} (q;q^2)_n (-aq^{2n+2};q)_\infty}{(q^2;q^2)_n}\\
  &=aq\sum_{n=0}^\infty 
   \frac{a^{2n+1} q^{2n^2+4n+1} (q;q^2)_{n} (1-q^{2n+1})(-aq^{2n+3};q)_\infty}{(q^2;q^2)_{n} }\\
  & \qquad\qquad +\sum_{n=0}^\infty 
   \frac{a^{2n+1} q^{2n^2+4n+1} (q;q^2)_n (-aq^{2n+2};q)_\infty}{(q^2;q^2)_n}\\
  &= \sum_{n=0}^\infty 
   \frac{a^{2n+1} q^{2n^2+4n+1} (q;q^2)_{n} 
(-aq^{2n+3};q)_\infty}{(q^2;q^2)_{n} } 
   \\ &\qquad\qquad\times     \Big( aq(1-q^{2n+1}) + (1+aq^{2n+2}) \Big)\\
  &=(1+aq) \sum_{n=0}^\infty 
   \frac{a^{2n+1} q^{2n^2+4n+1} (q;q^2)_{n} (-aq^{2n+3};q)_\infty}{(q^2;q^2)_{n} } \\
   &=aq(1+aq) F_1(aq),
\end{align*} and thus~\eqref{Feq3} is established.
\end{pf}

\begin{lem} \label{Elemma}The functions $E_1(a)$, $E_2(a)$, $E_3(a)$ satisfy the following
system of $q$-difference equations:
\begin{align}
  E_1(a) &= E_3(aq)\label{Eeq1}\\
  E_2(a) &= (1+aq) E_2(aq)\label{Eeq2}\\
  E_3(a) &= E_1(a) + aq(1+aq) E_1(aq) \label{Eeq3}
\end{align}
\end{lem}
\begin{pf}
Equations~\eqref{Eeq1} and~\eqref{Eeq2} are immediate from the definitions~\eqref{E1}--\eqref{E3}.
\begin{align*}
E_3(a) - E_1(a) &= \sum_{n=0}^\infty \sum_{m=0}^\infty \frac{ a^{2n+m} q^{2n^2+m^2+2mn}}{(q^2;q^2)_n (q;q)_m} (1-q^{2n+m})\\
&= \sum_{n=0}^\infty \sum_{m=0}^\infty \frac{ a^{2n+m} q^{2n^2+m^2+2mn}}{(q^2;q^2)_n (q;q)_m} 
  \Big( (1-q^{2n}) + q^{2n}(1-q^m) \Big)\\
&=\sum_{n=1}^\infty \sum_{m=0}^\infty \frac{ a^{2n+m} q^{2n^2+m^2+2mn}}{(q^2;q^2)_{n-1} (q;q)_m} 
+\sum_{n=0}^\infty \sum_{m=1}^\infty \frac{ a^{2n+m} q^{2n^2+m^2+2mn+2n}}{(q^2;q^2)_n (q;q)_{m-1}}\\
 &=\sum_{n=0}^\infty \sum_{m=0}^\infty \frac{ a^{2n+m+2} q^{2n^2+2n+m^2+2mn+2m+2}}{(q^2;q^2)_{n} (q;q)_m} \\ &\qquad\qquad
+\sum_{n=0}^\infty \sum_{m=0}^\infty \frac{ a^{2n+m+1} q^{2n^2+m^2+2m+2mn+4n+1}}{(q^2;q^2)_n (q;q)_{m}}\\
&=\sum_{n=0}^\infty \sum_{m=0}^\infty \frac{ a^{2n+m+1} q^{2n^2+m^2+2m+2mn+4n+1}}{(q^2;q^2)_n (q;q)_{m}} (1+aq)\\
&=aq(1+aq) \sum_{n=0}^\infty \sum_{m=0}^\infty \frac{ a^{2n+m} q^{2n^2+m^2+2m+2mn+4n}}{(q^2;q^2)_n (q;q)_{m}}\\
&=aq(1+aq)E_1(aq),
\end{align*} and thus~\eqref{Eeq3} is established.
\end{pf}

\begin{thm} 
\begin{equation}
 \sum_{n=0}^\infty T(n) q^n = 
    \sum_{n=0}^\infty \frac{q^{2n(n+1) } (q;q^2)_n (-q^{2n+ 2};q)_\infty}
 {(q^2;q^2)_n }
 \end{equation}
\end{thm}
\begin{pf}
  By a standard argument (see, e.g.~\cite[p. 442 ff., Lemma 1 and the remark following Eq. (5.9)]{gea:qdiff}), the system~\eqref{Feq1}--\eqref{Feq3} has unique solutions with $F_1(0)=F_2(0)=F_3(0)$.
  Since by Lemma~\ref{Elemma}, \eqref{E1}--\eqref{E3} satisfy the same system of $q$-difference
equations and $E_1(0)=E_2(0)=E_3(0)$, it follows that $E_i(a) = F_i(a)$, for $i=1,2,3$.
In particular, $E_1(1) = F_1(1)$.
  Further, observing that the inner sum in $E_1(1)$, by~\eqref{rr2}, generates partitions
  into parts which are distinct, nonconsecutive integers greater than $2n+1$,
  \[ \sum_{n=0}^\infty T(n) q^n = F_1(1), \]
  and thus the result follows.
\end{pf}

\begin{rem}
The motivation behind using the set of $q$-difference equations~\eqref{Feq1}--\eqref{Feq3}
comes from the fact that it is known (see~\cite{gea:qdiff},~\cite{dmb:even},~\cite{avs}) 
that the functions arising in $k=3$ case of Bressoud's theorem must satisfy~\eqref{Feq1}--\eqref{Feq3}.
Namely, it can be shown that if $Q_{i}(a) := (-aq;q)_\infty J_{1,i/2} ( ; a^2; q^2)$, where
$J_{k,i}(a;x;q)$ is defined in~\cite{gea:qdiff} and~\cite[p. 106, Eq. (7.2.2)]{top}, the $Q_{i}(a)$ for
$i=1,2,3$ also satisfy~\eqref{Feq1}--\eqref{Feq3}.
\end{rem}

\section{Discussion}
The $f$ map was created specifically to map the set $\Gtwo$ onto the set $\Atwo$.
As was shown, $f$ was also useful in a variety of other closely related contexts.  It might
be interesting to explore whether this map, or a generalization of it, might be applicable
in additional settings.

  This having been said, it is not clear what these additional settings
might be.
  Recall that our work began with the observation that
\begin{equation} \label{mod7}
\sum_{n=0}^\infty G_{3,2}(n) q^n = \underset{j\not\equiv 0,\pm 2{\hskip-3mm}\pmod{7}}
{\prod_{j\geqq 1}} \frac{1}{1-q^j} = 
\sum_{j=0}^\infty \frac{ q^{2j^2+2j} (-q^{2j+1};q)_\infty }{(q^2;q^2)_j},
\end{equation}
where the first equality follows from Gordon's theorem, and the
second equality is the second Rogers-Selberg identity.  Our $f$ map
provides a bijection between the partitions enumerated on the
left hand side with those enumerated in  
Andrews' combinatorial interpretation of the right.  The
middle member of~\eqref{mod7} serves only as a bridge between the
left and right and does not play a role in the bijection. 

 A natural modulus 5 analog of~\eqref{mod7} is
\begin{equation}
\label{mod5}
\sum_{n=0}^\infty G_{2,1}(n) q^n = \underset{j\not\equiv 0,\pm 1{\hskip-3mm}\pmod{5}}
{\prod_{j\geqq 1}} \frac{1}{1-q^j} =
\sum_{j=0}^\infty \frac{ (-1)^j q^{3j^2-2j} (-q^{2j+1};q)_\infty }{(q^2;q^2)_j},
\end{equation}
where the first equality follows from Gordon's theorem and the second equality
is an identity due to Rogers~\cite[p. 339, Ex. 2]{ljr1894} which also appears
in Slater's list~\cite[p. 154, Eq. (19)]{ljs}.
We are confronted with several difficulties immediately.  First,
the partitions enumerated by $G_{2,1}(n)$ do not have any repeated
parts and therefore our $f$ map becomes the identity map in this
context.  Next, even though the right hand side of~\eqref{mod5} resembles
that of~\eqref{mod7}, the presence of the factor $(-1)^j$ in~\eqref{mod5}
precludes the possibility of a combinatorial interpretation as simple
as that of the right hand side of~\eqref{mod7}.

  Perhaps the following modulus 9 analog of~\eqref{mod9} would be
more amenable to study via the methods of this paper:
\begin{multline}
\label{mod9}
\sum_{n=0}^\infty G_{4,1}(n) q^n = \underset{j\not\equiv 0,\pm 1{\hskip-3mm}\pmod{9}}
{\prod_{j\geqq 1}} \frac{1}{1-q^j} \\ =
\sum_{j=0}^\infty \frac{ q^{3j^2+3j} (-q^{3j+3};q^3)_\infty }{(q^3;q^3)_j 
  (1-q^{3j+2})} \prod_{h=j+1}^\infty \left( 1 + \frac{q^{3h+1}}{1-q^{3h+1}}
+ \frac{q^{3h+2}}{1-q^{3h+2}} \right),
\end{multline} 
where the first equality follows from Gordon's theorem and the second
equality is an identity due to Bailey~\cite[p. 422, Eq. (1.7)]{wnb} which
also appears in Slater~\cite[p. 156, Eq. (40)]{ljs}.
Andrews provides a (somewhat complicated) combinatorial interpretation
of the last member of~\eqref{mod9} in~\cite[\S 5]{gea}.
It should be noted, however, that in the $k=4$ case of Gordon's theorem, parts
may be repeated up to three times (with certain additional restrictions)
and thus it is not obvious how to adapt or generalize the $f$ map to  
this situation.
 
\section*{Acknowledgments}
The author thanks George Andrews for a number of stimulating and helpful 
discussions, and the referees for useful comments.

 \end{document}